\documentclass{amsart}   
 \usepackage{amsmath}
 \usepackage{amssymb}
 \usepackage{amscd}
 \usepackage[mathscr]{eucal}
 \newtheorem{thm}{Theorem}[section]

 \newtheorem{cor}[thm]{Corollary}
 \newtheorem{prop}[thm]{Proposition}
 \newtheorem{defn}[thm]{Definition}

 \newcommand{\Q}{\mathbb{Q}}

 \title{ $g$-elements of matroid complexes}  
 \author{Edward Swartz} 
 \address{Malott Hall \\ Cornell University \\ Ithaca, NY14853} 
 \email{ebs@math.cornell.edu}
 \thanks{Partially supported by a VIGRE postdoc under NSF grant number
9983660
to Cornell University.}

\begin{document}
 \begin{abstract}
  A $g$-element for a graded $R$-module is a one-form with properties 
  similar to 
  a Lefschetz class in the cohomology ring of a compact complex projective 
  manifold, except that the induced multiplication maps are injections 
  instead of bijections.  We show that if $k(\Delta)$ is the face ring of 
  the independence complex of a matroid and the characteristic of $k$ is 
  zero, then there is a non-empty Zariski open subset of pairs 
  $(\Theta,\omega)$ such that $\Theta$ is a linear set of paramenters  
  for $k(\Delta)$ and $\omega$ is a $g$-element for 
  $k(\Delta)/<\Theta>.$    This leads to 
  an inequality on the first half of the $h$-vector of the complex similar 
  to the $g$-theorem for simplicial polytopes.

 \end{abstract}

 \maketitle
\section{Introduction}
   The combinatorics of the independence complex of a matroid can be 
   approached from several different directions.  The $f$-vector directly 
   encodes the number of independent sets of every cardinality, while the 
   $h$-vector contains the same information encoded in a way which is 
   more appropriate for reliability problems \cite{CC}.  In either case 
   the fundamental question is the same.  What vectors are possible?  

  Let $(h_0, h_1, \dots h_r)$ be the $h$-vector of the independence 
  complex of a rank $r$ matroid without coloops.  Using a PS-ear 
  decomposition of the complex Chari \cite{Ch} proved that for all 
  $i \le r/2, \ h_i 
  \le h_{r-i}$ and $ h_{i-1} \le h_i$ . By showing that the 
  $h$-vector was the Hilbert function of $k(\Delta)/<\Theta>$, where 
  $k(\Delta)$ is the face ring of the complex and $\Theta$ is a linear 
  system of parameters for $k(\Delta),$ 
  Stanley \cite{St1} proved that $h_{i+1} \le h_i^{<i>}$ (see section 
  \ref{face 
  rings} for a definition of the $<i>$-operator). By combining 
  these two methods we show in Theorem \ref{combinatorics} that if we 
  define $g_i = 
  h_i - h_{i-1},$ then $g_{i+1} \le g_i^{<i>}$ for all $i < r/2.$  
  All 
  of these inequalities are immediate consequences of the existence of 
  pairs $(\Theta,\omega)$ such that $\Theta$ is a linear set of parameters 
  for $k(\Delta)$ and $\omega$ is a $g$-element for 
  $k(\Delta)/<\Theta>.$  Using a different approach, toric hyperk\"{a}hler 
  varieties, Hausel and Sturmfels proved the existience of $g$-elements 
  for $k(\Delta)/<\Theta>$ when the matroid is representable over the 
  rationals \cite{HS}.
   A $g$-element is a one-form which acts like a 
  Lefschetz class of a compact complex projective manifold except that it 
  induces injections instead of bijections (Definition \ref{g element}).

  The broken circuit complex of a matroid is subcomplex of the 
  independence complex and directly encodes the coefficients of the 
  characteristic polynomial of the matroid .  Every broken circuit 
  complex is a cone, and if we remove the cone point we obtain a reduced 
  broken circuit complex.  Any independence complex is also a reduced broken 
  circuit complex.  Since the $h$-vector is unchanged by the removal of a 
  cone point, the set of $h$-vectors of independence complexes is a 
  (strict) subset of the set of $h$-vectors of broken circuit complexes.  
  A natural question is whether or not Theorem \ref{main} holds for broken 
  circuit complexes.  In Section \ref{BC complex} we show that even if 
  broken circuit complexes satisfy the corresponding combinatorial 
  inequalities, there may be no set of linear parameters for the face ring 
  such 
  that there exist $g$-elements for the quotient ring.

  Matroid terminology and notation closely follows \cite{O}.   The main 
  exception to this is that we use $M-A$ for the deletion of a subset 
  instead of $M \setminus A.$ The ground set of the matroid $M$ is always 
  $E.$ 
 
 \section{Complexes}

Let $\Delta$ be a finite abstract simplicial complex with vertices $V = 
\{v_1,\dots,v_n\}.$  The {\it f-vector} of $\Delta$ is the sequence 
$(f_0(\Delta),\dots,f_s(\Delta)),$ where $f_i(\Delta)$ is the number 
of simplices of cardinality 
$i$ and $s-1$ is the dimension of $\Delta.$  The {\it h-vector} of 
$\Delta$ is the sequence $(h_0(\Delta),\dots,h_s(\Delta))$ defined by,
$$ h_i(\Delta) = \sum^i_{k=0} (-1)^{i+k} f_k(\Delta) \binom{s-k}{i-k}.$$
Equivalently,  if we let $f_\Delta(t) = f_0 t^s + f_1 t^{s-1} + \dots 
f_{s-1} t + f_s,$ then $h_\Delta(t) = h_0 t^s + h_1 t^{s-1} + \dots 
h_{s-1} t + h_s$ satisfies $h_\Delta(1+t) = f_\Delta(t).$

The {\it independence complex} of $M$ is the simplicial complex whose 
vertices are the non-loop elements of $E$ and whose simplices are the 
independent subsets of $E.$  We let $\Delta(M)$ represent the independence 
complex of $M.$

In order to define the broken 
circuit complex for $M$,
we first choose a linear order ${\bf n}$ on the elements of the matroid.  
Given such 
an order, a {\it broken circuit} is a circuit with its least element 
removed. 
The {\it broken circuit complex} is the simplicial complex whose simplices 
are the subsets of $E$ which do not contain a broken circuit.  
We denote the broken circuit complex of $M$ and ${\bf n}$ by 
$\Delta^{BC}(M),$ or $\Delta^{BC}(M,{\bf n}).$  
Different orderings may lead to different complexes, see \cite[Example 
7.4.4]{Bj}. Conversely, distinct matroids can have the same broken circuit 
complex.  For instance, let $E=\{e_1,e_2,e_3,e_4,e_5,e_6\},$ and let 
${\bf n}$ be the obvious order.  Let $M_1$ be the matroid on $E$ whose 
bases are all triples except $\{e_1,e_2,e_3\}$ 
and  $\{e_4,e_5,e_6\}$ and let $M_2$ be the matroid on $E$ whose 
bases are all triples except $\{e_1,e_2,e_3\}$ and 
$\{e_1,e_5,e_6\}.$  Then $M_1$ and $M_2$ are non-isomorphic matroids but 
their broken circuit complexes are identical.  Both $h_{\Delta(M)}(t)$ and 
$h_{\Delta^{BC}(M)}(t)$ satisfy similar contraction-deletion formulas.

\begin{prop} \cite{Bj}, \cite{BO} \label{tuttepoly}
 \begin{enumerate}
  \item
    If $e$ is a loop of $M,$ then $h_{\Delta(M)}(t) = 
    h_{\Delta(M-e)}(t),$ and $h_{\Delta^{BC}(M)}(t)=1.$ 
  \item
    If $e$ is a coloop of $M,$ then $h_{\Delta(M)}(t) = 
    h_{\Delta(M-e)}(t),$ and $h_{\Delta^{BC}(M)}(t)=h_{\Delta^{BC}(M-e)}(t).$
  \item
    If $e$ is neither a loop nor a coloop of $M,$ then
    $h_{\Delta(M)}(t) = h_{\Delta(M-e)}(t) + h_{\Delta(M/e)}(t)$ and
    $h_{\Delta^{BC}(M)}(t) = h_{\Delta^{BC}(M-e)}(t) + 
    h_{\Delta^{BC}(M/e)}(t).$ 
   \item
    If $M = M_1 \oplus M_2,$ then $h_{\Delta(M)}(t) = h_{\Delta(M_1)}(t) 
    \cdot 
     h_{\Delta(M_2)}(t).$
   \item
If $S$ is a series class of $M$ which is not a circuit, 
then $h_{\Delta(M)}(t) = h_{\Delta(M/S)}(t) + h_{\Delta(M-S)}(t) 
(1+t+\dots+t^{|S|-1}).$
\end{enumerate}
\end{prop}
 
 \section{Face rings} \label{face rings}

Let $k$ be a field and let $R=k[x_1,\dots,x_n].$

\begin{defn}
  The {\it face ring} of $\Delta$ is the graded $k$-algebra
$$k[\Delta] = R/I_\Delta,$$
where $I_\Delta$ is the ideal generated by all monomials $x_{i_1} 
\cdot \dots \cdot x_{i_l}$ such that $\{v_{i_1},\dots,v_{i_l}\}$ is not a 
face of $\Delta.$  
\end{defn}

Let $s-1$ be the dimension of $\Delta.$ 
Let $\Theta = \{\theta_1,\dots,\theta_s\}$ be a set of one-forms in $R.$  
Write each $\theta_i = k_{i1} x_1 + \dots k_{is} x_s$ and let $K = (k_{ij}).$  
To each simplex in $\Delta$ there is a corresponding set of column vectors 
in $K.$  If for every simplex of $\Delta$ the corresponding set of column 
vectors is independent, then $\Theta$ is a {\it linear set of parameters} 
(l.s.o.p.) 
for $k(\Delta).$  If $k$ is infinite, then it is always possible to choose $\Theta$ 
such that every set 
of $s$ columns of $K$ is independent. 

Given a l.s.o.p. $\Theta$ for $k(\Delta)$ let $R(\Delta,\Theta) = 
k(\Delta)/<\Theta>.$ If $\Theta$ is unambiguous, then we just use 
$R(\Delta).$ Since $\Theta$ is homogeneous $R(\Delta)$ is a 
graded $k$-algebra.

\begin{thm} \cite{St1}
  Let $\Theta$ be a l.s.o.p. for $\Delta(M)$ and let $R(\Delta(M))_i$ be 
  the $i^{th}$ graded component of $R(\Delta(M)).$  Then $h_i(\Delta(M)) = 
  \dim_k R(\Delta(M))_i.$  Similarly, if $\Theta$ is a l.s.o.p. for 
  $\Delta^{BC}(M),$ then $h_i(\Delta^{BC}(M)) = \dim_k R(\Delta^{BC}(M))_i.$
\end{thm}

Given any two integers $i,j > 0$ there is a unique way to write 
$$j=\binom{a_i}{i} + \binom{a_{i-1}}{i-1} + \dots + \binom{a_l}{l}, a_i > 
a_{i-1} > \dots > a_l \ge l \ge 1.$$
Given this expansion define,
$$j^{<i>} = \binom{a_i+1}{i+1} + \binom{a_{i-1}+1}{i} + \dots + \binom{a_l+1}{l+1}, a_i > 
a_{i-1} > \dots > a_l \ge l \ge 1.$$

\begin{thm} \cite[Theorem 2.2]{St} \label{stanley}
  Let $Q=R/I,$ where $I$ is a homogeneous ideal.  Let $Q_i$ be the forms of degree $i$ in $Q$ and let $h_i = \dim_k Q_i.$  
  Then $h_{i+1} \le h_i^{<i>}.$
\end{thm}

\begin{cor} \cite{St1}
  For any independence or broken circuit complex $h_{i+1} \le h_i^{<i>}.$
\end{cor}

\section{The ring $R(\Delta(M))$}

In order to study the properties of $h_i(\Delta(M))$ we will look for 
elements with properties slightly weaker than those provided by Lefschetz 
elements of the cohomology ring of a compact complex projective manifold. 

\begin{defn} \label{g element}
  Let $N$ be a (non-negatively) graded $R$-module whose dimension over $k$ 
  is finite.   
  Let $r$ be the last non-zero grade of $N$ and let $\omega$ be a one-form in $R.$  
  Then $\omega$ is a {\it $g$-element} for $N$ if for all $i, 0 \le i \le 
  r/2,$ multiplication by $\omega^{r-2i}$ is an injection from $N_i$ to $N_{r-i}.$
\end{defn}

\noindent  If we replace injection with bijection in the above definition, 
then we obtain the strong Stanley property in \cite{Wa}.

Let $M$ be a rank $r$ matroid without coloops and $k$ a field of characteristic 
  zero.  Let $n=|E|.$  Write the elements of 
  $k^{n \times (r+1)}$ in the form $(\Theta, \omega),$ 
  where $\Theta$ consists of $r$ elements in $k^n$ and $\omega$ is also in 
  $k^n.$ Identify elements of 
  $k^n$ with the one-forms in $R$ in the canonical way.  Let  
  $U$ be the set of all pairs $(\Theta,\omega) \in k^{n \times (r+1)}$ 
  such that $\Theta$ is a l.s.o.p. 
  for $k(\Delta(M))$ and $\omega$ is a $g$-element for 
  $R(\Delta(M),\Theta).$
\begin{thm} \label{main}
 Let $M,U$ and $k$ be as above.  Then, $U$ is 
 a non-empty Zariski open subset of $k^{n \times (r+1)}.$
\end{thm}

 \begin{proof}
  We first note that $\Theta$ is a l.s.o.p. for $k(\Delta)$ if and only if 
  the determinants of the appropriate $r \times r$ minors of $K$ are 
  non-zero.  Secondly, the multiplication maps $\omega^{r-2i}$ 
  can be encoded as matrices 
  which are polynomial in the coefficients of $K$ and $\omega.$  Thus, $U$ 
  is  the intersection of two Zariski open subsets of $k^{n \times (r+1)}.$
 
  To show that $U$ is not empty we proceed by induction on $n.$  However, 
  we use a slightly different (but equivalent) induction hypothesis. Let 
  $C(j)$ be the circuit with $j$ elements.  Let $P$ be a direct sum of 
  circuits, so we can write $P = C(j_1) \oplus \dots \oplus C(j_m).$  The 
  rank of $M \oplus P$ is $r^\prime=r+j_1+ \dots +j_m -m$ and its 
  cardinality is $n^\prime = n + j_1 + \dots + j_m.$ The induction 
  hypothesis is that given any such $P,$
  then $U = \{(\Theta,\omega) \in k^{n^\prime \times (r^\prime +1)}: 
  \Theta$ is a l.s.o.p. for $k(\Delta(M \oplus P))$ and 
  $\omega$ is a $g$-element for $R(\Delta( M \oplus P),
  \Theta) \}$ is not empty.  If $M$ consists of a single 
  loop, then $k(\Delta(M \oplus P)) \simeq k(\Delta(P)).$ 
  As a simplicial complex $\Delta(P)$ is $\partial (\Delta^{j_1-1}) 
  \ast \dots \ast \partial (\Delta^{j_m-1}),$ where $\Delta^j$ is the 
  $j$-simplex.  Since this is the boundary of a convex rational 
  polytope we can 
  apply the Hard Lefschetz theorem as in \cite{St2} to see that $U$ is 
  not empty when $k= \Q,$  and hence is not empty for any field of 
  characteristic zero. 

Suppose that $M$ consists of a single coloop.  Given $P$ and 
$\Theta$ as in the above paragraph, $\Theta \cup \{x_{n^\prime}\}$ 
is a l.s.o.p. 
for $k(\Delta(P \oplus M))$ and $\omega,$ which can be 
viewed as an 
element of $k[x_1,\dots,x_{n^\prime}],$ is a $g$-element for $R(\Delta(P 
\oplus M),\Theta \cup \{x_{n^\prime}\}).$  

For the induction step, let $S$ be a series class of $M.$  Reordering $M$ 
if necessary, we assume that $S=\{e_1,\dots,e_s\}$ consists of the first 
$s$ elements of $M.$
If $S$ is a 
circuit, then $M = (M-S) \oplus C(s).$ Hence, $M \oplus P = (M-S) \oplus (C(s) 
\oplus P)$ and the 
induction hypothesis applies to $M-S.$  So assume that $S$ is 
independent .   Let $x^S = x_1 \cdot \dots \cdot x_s.$  For 
$\Theta$ a l.s.o.p. for $k(\Delta(M \oplus P))$ consider the following 
short exact sequence.

\begin{equation} \label{ses}
 0 \to \frac{<x^S> R}{<x^S> \cap (\Theta + I_{\Delta(M \oplus P)})} \to 
R(\Delta(M \oplus P,\Theta)) \to 
\frac{R(\Delta( M \oplus P),\Theta)}{<x^S>} \to 0.
\end{equation}

Since $S$ is a series class, a subset of $M-S$ is independent if and only 
if its union with any proper subset of $S$ is independent.  Hence, the 
right-hand side is just $R(\Delta((M-S) \oplus (C(s) \oplus 
P),\Theta).$ Therefore,  we can apply the induction hypothesis to $M-S$ to 
obtain a non-empty 
Zariski open subset $U^\prime$ of $k^{n \times (r+1)}$ consisting 
of pairs $(\Theta^\prime, \omega^\prime)$ such that $\omega^\prime$ is 
a $g$-element for $R(\Delta((M-S) \oplus C(s) \oplus P), \Theta^\prime).$

In order to analyze the left-hand side of (\ref{ses}) choose generators 
$\{\theta_1, \dots, \theta_s, \dots, \theta_{r^\prime}\}$ for $<\Theta>$ 
so that in the corresponding matrix $K, k_{ij} = \delta_{ij}$ for $1 \le 
i \le s.$  
Now define an $R$-module structure on $R^\prime = k[x_{s+1},\dots,x_{n^\prime}]$ 
by defining $(x_i) \cdot f = (x_i - \theta_i) \cdot f$ for 
$1 \le i \le s$ and $f \in R^\prime.$   
Let $\phi: R^\prime \to 
(<x^S> R)/(\Theta + I_{\Delta(M \oplus P)})$ be multiplication 
by $x^S.$  Since $S$ is independent, every polynomial in 
$(<x^S> R)/(\Theta + I_{\Delta(M \oplus P)})$ is equivalent to 
a polynomial in $<x^S> R^\prime.$  So, $\phi$ is 
surjective.   
The kernel of $\phi$ contains
$\Theta^\prime = \{\theta \in \Theta: 
\theta = k_{s+1} x_{s+1} + \dots + k_{n^\prime} x_{n^\prime}.\}$ 
In addition, $\ker \phi$ contains all monomials in 
$I_{\Delta((M/S) \oplus P)}.$ Since  $\Theta^\prime$ is a l.s.o.p. 
for $k(\Delta(M/S)),$ we see that $\phi$ is a degree $s$ graded
surjective $R$-module 
homomorphism from $R^\prime/(I_{\Delta(( M/S) \oplus P)} + 
\Theta^\prime)$ to the left-hand side of (\ref{ses}) .   
Proposition \ref{tuttepoly} and  $h_{\Delta(C(s))}(t) = 
1+t+\dots+t^{s-1}$ show that the $k$-dimension of 
$R^\prime(\Delta(M/S),\Theta^\prime)$ and the l.h.s. of (\ref{ses}) are the 
same.  Hence $\phi$ is an isomorphism. Therefore, by the induction 
hypothesis applied to $M/S,$ there is a non-empty Zariski open subset 
$U^{\prime \prime}$ of $k^{n \times (r+1)}$ consisting of pairs 
$(\Theta^{\prime \prime},\omega^{\prime \prime})$ such that the 
multiplication map

$$\omega^{\prime \prime (r^\prime - 2i -s)}: 
(\frac{<x^S> R}{<x^S> \cap (\Theta^{\prime 
\prime} + I_{\Delta(M \oplus P)})})_{i+s} \to 
(\frac{<x^S> R}{<x^S> \cap (\Theta^{\prime 
\prime} + I_{\Delta(M \oplus P)})})_{r^\prime - i}$$
is an injection for $1 \le i \le (r^\prime -s)/2.$  Now, $U^\prime \cap 
U^{\prime \prime} \subseteq U.$  Since the intersection of two non-empty 
Zariski open 
subsets of $k^{n \times (r+1)}$ is not empty, $U$ is also not empty. 

\end{proof}
  
\begin{thm} \label{combinatorics}
Let $M$ be a rank $r$ matroid without coloops.  Let $h_i = 
h_i(\Delta(M)).$  Then,
\begin{enumerate}
  \item
     $h_0 \le \dots \le h_{\lfloor r/2 \rfloor}.$

  \item
     $h_i \le h_{r-i}$ for all $i \le r/2.$

   \item
     Let $g_i = h_i - h_{i-1}.$  Then, for all $i < r/2,\  g_{i+1} \le 
     g_i^{<i>}.$
\end{enumerate}
\end{thm}

 \begin{proof}
  The first two inequalities follow from the injectivity properties of any 
  $g$-element $\omega$ for $R(\Delta(M),\Theta).$  Since $g_i = 
  (R(\Delta(M),\Theta)/<\omega>)_i$ when $i < r/2,$ the last inequality 
  follows from Theorem \ref{stanley}.
\end{proof}

The first two inequalities were obtained by Chari using a PS-ear 
decomposition of $\Delta(M).$  See \cite{Ch2} for details on PS-ear 
decompositions. Hausel and Sturmfels used toric hyperk\"{a}hler varities to 
prove the last inequality for matroids representable over the rationals 
\cite{HS}.
The proof of Theorem \ref{main} is essentially an algebraic 
version of a PS-ear decomposition \cite[Theorem 2]{Ch}. Indeed, the 
proof works without change for any simplicial complex with a PS-ear 
decomposition.


\section{The ring $R(\Delta^{BC}(M))$} \label{BC complex}

As shown in \cite{Bry4} the cone on any independence complex is a broken 
circuit complex (for some other matroid).  Since the $h$-vector of the 
cone of a simplicial complex is the same as the $h$-vector of the 
original complex, the $h$-vectors of independence complexes form a 
(strict) subset of the the $h$-vectors of broken-circuit complexes.  It 
is natural to ask whether or not Theorem \ref{main} holds for 
$\Delta^{BC}(M).$  The last non-zero element of the $h$-vector of 
$\Delta^{BC}(M)$ is $r-m,$ 
where $m$ is the number of components of $M.$  It is not difficult to 
modify the proof of Theorem \ref{main} to produce injections from 
$R(\Delta^{BC}(M))_0$ to $R(\Delta^{BC}(M))_{r-m}$ and from $R(\Delta^{BC}(M))_1$ 
to $R(\Delta^{BC}(M))_{r-m-1}.$ Since the first possible problem is in degree 
2, the smallest possible rank of $M$ for which Theorem \ref{main} does 
not hold  for $\Delta^{BC}(M)$ is six.

Let $G(s)$ be the graph obtained by subdividing each edge of the graph 
consisting of $s$ parallel edges into two edges.  Let $M(s)$ be the cycle 
matroid of $G(s).$  The rank of $M(s)$ is $s+1$ and $M(s)$ has $2s$ 
elements.

\begin{prop}
  Let $\Theta$ be a l.s.o.p. for $M(s), {\bf n}$  a linear order of 
  the elements of $M(s)$  and $\omega$ a linear form in $k[x_1,\dots, 
  x_{2s}].$ Then, multiplication by $\omega$ has 
  a non-trivial kernel in $R(\Delta^{BC}(M))_2.$
\end{prop}

\begin{proof}
  Let $E_l$ consist of the greatest $l$ elements of $M(s)$ with respect to 
  ${\bf n}.$  Let $\{e_i,e_j\}$ be the first pair of edges to appear in 
  $E_l$ as $l$ goes from $1$ to $s+1$ such that they come from the 
  subdivision of one of the parallel edges used to construct $G(s).$  
  Consider the ideal $<x_i x_j> \subseteq R(\Delta^{BC}(M,{\bf n})).$  
  Using the same reasoning as in the proof of Theorem \ref{main}, the 
  choice of $\{e_i,e_j\}$ implies that $<x_i x_j>$ is isomorphic as an 
  $R$-module to $R^\prime(\Delta^{BC}(\Delta(M(s)/\{e_i,e_j\},{\bf 
  n}^\prime),\Theta^\prime),$ 
  where $R^\prime$ and $\Theta^\prime$ are defined as in the proof of 
  Theorem \ref{main}, and ${\bf n}^\prime$ is the order on 
  $M(s)/\{e_i,e_j\}$ induced from ${\bf n}.$   Now, $M(s)/\{e_i,e_j\}$ is 
  the cycle matroid of the $G(s)$ with the two edges $\{e_i,e_j\}$ 
  contracted.  For any such pair and any linear order 
  $\Delta^{BC}(M(s)/\{e_i,e_j\}$ is an $s-2$ dimensional simplex.  Hence 
  $<x_i x_j> \simeq k$ and will vanish under any multiplication map.  
\end{proof}

Repeated application of Proposition \ref{tuttepoly} shows that 
$h_i(\Delta^{BC}(M(s))) = \binom{s}{i}$ when $i \neq 1$ and 
$h_1(\Delta^{BC}(M(s))) = s-1.$  When $s \ge 5,$  the $h$-vector 
of the broken circuit complex of $M(s)$ satisfies the combinatorial conditions of Theorem 
\ref{combinatorics} but there is no l.s.o.p. for the face ring such that 
the quotient ring has $g$-elements.  As far as we know, whether or not 
broken circuit complexes satisfy the combinatorial inequalities of 
Theorem \ref{combinatorics} remains an open question.
\\
\noindent {\it Acknowledgement:}  Louis Billera suggested the problem of 
determining whether or not independence complexes satisfied $g_{i+1} \le 
g_i^{<i>}.$

\end{document}